\documentclass[graybox]{svmult}

\usepackage{mathptmx}       
\usepackage{helvet}         
\usepackage{courier}        
\usepackage{type1cm}        
\usepackage{makeidx}         
\usepackage{graphicx}        
\usepackage{multicol}        
\usepackage[bottom]{footmisc}
\usepackage{amsfonts}
\usepackage{amsmath}
\usepackage{amsopn} 
\usepackage[hardware,software]{mymacros}

\usepackage{cleveref}
\usepackage{pgfplots}

\usetikzlibrary{shapes,snakes}
\usetikzlibrary{arrows}
\usepackage{tikz}
\usepackage{pgfplots}

\usepackage[utf8]{inputenc}
\usepackage[overload]{empheq}
 \usepackage{cases} 
 
\newenvironment{customlegend}[1][]{%
	\begingroup
	\csname pgfplots@init@cleared@structures\endcsname
	\pgfplotsset{#1}%
}{%
\csname pgfplots@createlegend\endcsname
\endgroup
}%
\def\addlegendimage{\csname pgfplots@addlegendimage\endcsname}

\usepackage[ruled,linesnumbered]{algorithm2e}

\makeindex             

\begin{document}

\title*{An Iterative Model Reduction Scheme for Quadratic-Bilinear Descriptor Systems with an Application to Navier-Stokes Equations}
  \titlerunning{A MOR Scheme for QB Systems with an Application to NS Equations}
\author{Peter Benner and Pawan Goyal}
\institute{Peter Benner \at Max Planck Institute for Dynamics of Complex Technical Systems, Sandtorstra\ss e 1, 39106, Magdeburg, Germany, \email{benner@mpi-magdeburg.mpg.de}
\and Pawan Goyal \at Max Planck Institute for Dynamics of Complex Technical Systems, Sandtorstra\ss e 1, 39106, Magdeburg, Germany, \email{goyalp@mpi-magdeburg.mpg.de}
}
%
%
\maketitle

\abstract*{We discuss model reduction for a particular class of quadratic-bilinear (QB) descriptor systems.  The main goal of this article is to extend the recently studied interpolation-based optimal model reduction  framework for QBODEs [Benner et al. '16] to a  class of descriptor systems in an efficient and reliable way. Recently, it has been shown in the case of linear or bilinear systems that   a direct extension of interpolation-based model reduction techniques to descriptor systems, without any modifications, may lead to poor reduced-order systems.  Therefore, for the analysis, we aim at transforming the considered QB descriptor system into an equivalent QBODE system by means of projectors for which  standard model reduction techniques for QBODEs can be employed, including aforementioned interpolation scheme. Subsequently, we discuss related computational issues,   thus resulting in a modified algorithm that allows us to construct \emph{near}--optimal reduced-order systems without explicitly computing the projectors used in the analysis.	The efficiency of the proposed algorithm is illustrated by means of a numerical example, obtained via  semi-discretization of the Navier-Stokes equations.}
\abstract{We discuss model reduction for a particular class of quadratic-bilinear (QB) descriptor systems.  The main goal of this article is to extend the recently studied interpolation-based optimal model reduction  framework for QBODEs [Benner et al. '16] to a  class of descriptor systems in an efficient and reliable way. Recently, it has been shown in the case of linear or bilinear systems that   a direct extension of interpolation-based model reduction techniques to descriptor systems, without any modifications, may lead to poor reduced-order systems.  Therefore, for the analysis, we aim at transforming the considered QB descriptor system into an equivalent QBODE system by means of projectors for which  standard model reduction techniques for QBODEs can be employed, including aforementioned interpolation scheme. Subsequently, we discuss related computational issues,   thus resulting in a modified algorithm that allows us to construct \emph{near}--optimal reduced-order systems without explicitly computing the projectors used in the analysis.	The efficiency of the proposed algorithm is illustrated by means of a numerical example, obtained via  semi-discretization of the Navier-Stokes equations.}

\section{Introduction}
\label{sec:intro} 
We investigate model order reduction for quadratic-bilinear (QB) descriptor systemsxx of the form
\begin{subequations}\label{eq:NS_full}
\begin{align}
E_{11}\dot{v}(t) &= A_{11}v(t) + A_{12}p(t) + Hv(t)\otimes v(t) + \sum_{k=1}^mN_kv(t)u_k(t) + B_1u(t),\\
0 &= A_{21}v(t) + B_2u(t),\quad v(0) = v_0,\\
y(t) &= C_1v(t) + C_2p(t),
\end{align}
\end{subequations}
where $E_{11}$, $A_{11}$, $N_k \in \R^{n_v\times n_v}$, $H \in \R^{n_v\times n_v^2}$,   $A_{12}$, $A_{21}^T\in \R^{n_v\times n_p}$, $B_1\in \R^{n_v\times m}$, $B_2\in \R^{n_p\times m}$, $C_1\in \R^{p\times n_v}$, $C_2\in \R^{p\times n_p}$; $v(t) \in \R^{n_v}$ and $p(t)\in \R^{n_p}$ are the state vectors; $u(t)\in \Rm$ and $y(t)\in \Rp$ are the control input and measured output vectors of the system, respectively, and $u_k(t)$ is the $k$th component of the input vector; $v_0 \in \R^{n_v}$ is an initial value for $v(t)$. Furthermore, we assume that $E_{11}$ and $A_{21}E_{11}^{-1}A_{12}$ are invertible. Hence, the linear part of the system~\eqref{eq:NS_full} ($H =0$, $N_k = 0$) has an index-2 structure, e.g., see \cite{HaiW02}. The special structure for QB descriptor systems as in~\eqref{eq:NS_full}, particularly, appears in the modeling of the incompressible Navier-Strokes equations. 

A high-fidelity modeling of dynamical systems is required to have a better understanding of underlying dynamical behaviors of a system. However, numerical simulations of such high-fidelity systems are expensive and often inefficient.  Thus, it is not a straightforward task, or sometimes not even possible,  to perform engineering design studies using these high-fidelity systems. One approach to circumvent this problem is \emph{model order reduction} (MOR), aiming at constructing surrogate models (reduced-order models) which are less complicated and replicate the important dynamics of the high-fidelity system. 

MOR techniques for linear systems are now very well-established and are widely applied in numerous applications, e.g., see \cite{morAnt05,morBauBF14,morBenMS05,morSchVR08}. Several of those techniques have been successfully extended to special classes of nonlinear ODE systems, namely bilinear and QB systems, see, e.g.,~\cite{morBenB12b,morBenB12,morBenB15,morBenD11,morBTQBgoyal,flagg2015multipoint,morGu09,morZhaL02}. These techniques can be classified mainly into two categories, namely, trajectory-based methods, and system theoretic approaches such as balanced truncation and moment-matching.  Briefly mentioning the primary ideas of these methods,  trajectory-based methods rely on a set of snapshots of the state solutions for training inputs, which is  then  used to determine a \emph{Galerkin} projection, for more details, see, e.g.,~\cite{morAstWWetal08,morChaS10,morHinV05,morKunV02,morKunV08}. Furthermore, the idea of balanced truncation is to find states which are hard to control as well as hard to observe, and truncating such states gives us a reduced-order system. This method for QB systems has been  recently studied in \cite{morBTQBgoyal}.   Interpolation-based methods aim at constructing reduced-order systems which approximate the input-output behavior of the system.  With this intent, such a  problem for QB systems was first considered in \cite{morGu09}, where  a one-sided projection method  to obtain an interpolating reduced-order system is proposed.  Later on, a similar problem  was addressed  in \cite{morBenB12a,morBenB15}  for single-input single-output (SISO) QB systems, where a two-sided projection method was proposed, ensuring a higher number of moment to be matched. However, the main challenges for this method are a good selection of interpolation points and the application to multi-input multi-output (MIMO)  QB systems.  To address these issues, an interpolation-based optimal model reduction problem for QB systems was recently addressed in \cite{morBenGG16}, where a reduced-order system is constructed, aiming at minimizing a system norm of the error system, e.g., a truncated $\cH_2$-norm. 

However, there has been little attention  paid to descriptor systems (DAEs) which involve algebraic constraints as well along with a differential part, and this still requires further research. In the direction of MOR for nonlinear DAEs, interpolation-based methods for specially structured bilinear DAEs have been investigated, e.g., in \cite{morBenG16,goyal2016iterative}. Furthermore, a moment-matching method for SISO QB systems as in~\eqref{eq:NS_full} was studied in \cite{MPIMD15-18}. However, the first  challenge in this method is the choice of interpolation points which plays a crucial role in determining the quality of a reduced-order system, and secondly, it is applicable only to SISO systems which is certainly  very restrictive  from a real-world application point of view.   

 In this work, we aim to extend an $\cH_2$-optimal model reduction framework for QBODEs~\cite{morBenGG16} to QBDAEs, having the structure as in \eqref{eq:NS_full}. To that end, we first recall  an interpolation-based $\cH_2$-optimal model reduction technique for QBODEs  and the corresponding iterative scheme to construct reduced-order systems, in the subsequent section. In Section~\ref{sec:transferDAEtoODE}, we first present the transformation of the system~\eqref{eq:NS_full} into an equivalent ODE system by means of projectors. We further investigate how the iterative scheme can efficiently be applied to the equivalent ODE systems to obtain reduced-order systems without computing the projectors explicitly.  Finally, in Section~\ref{sec:numerics}, we illustrate the proposed methodology using a lid driven cavity model, which is obtained by semi-discretized Navier-Stokes equations.

\section{Model Reduction for Quadratic-Bilinear  ODEs and  Related Work}
\label{sec:BG_work}
In this section, we briefly discuss an $\cH_2$-optimal model reduction problem for QBODEs. We begin by introducing the problem setting for QBODEs. For this,
 we consider a QB system in the state-space form:
 \begin{subequations}\label{eq:QBODE_full}
\begin{align}[left={{\Sigma}:\empheqlbrace\,}]
\dot{x}(t) &= Ax(t) + Hx(t)\otimes x(t) + \sum_{k=1}^m N_k x(t) u_k(t) + Bu(t),\\
y(t) &  = Cx(t), \quad x(0)=0,
\end{align}
\end{subequations}
where $x(t) \in \Rn$, $u(t) \in \Rm$ and $y(t) \in \Rp$ are the state, input, and output vectors at the time instant $t$; all other matrices are real and are of appropriate dimensions. It is assumed that the matrix $A$ is stable. Furthermore,  one may also consider a  general nonsingular matrix $E$ in front of $\dot{x}(t)$; however, to keep the discussion simple in this section, we consider it to be an identity matrix.  

In the context of model reduction, our aim is to replace the system \eqref{eq:QBODE_full}  with a simpler and reliable reduced system, having a form:
 \begin{subequations}\label{eq:QBODE_red}
 	\begin{align}[left={\hat{\Sigma}:\empheqlbrace\,}]
 \dot{\hx}(t) &= \hA\hx(t) + \hH\hx(t)\otimes \hx(t) + \sum_{k=1}^m \hN_k \hx(t) u_k(t) + \hB u(t),\\
 	\hy(t) &  = \hC\hx(t), \quad \hx(t) = 0,
 	\end{align}
 \end{subequations}
   where $\hx(t)\in \Rr$ with $r\ll n$, all other reduced matrices are of appropriate dimensions, and the reduced output $\hy(t)$ approximates the corresponding original output $y(t)$ for a wide range of admissible $L_2$-bounded system inputs.   
 
 We focus on constructing a reduced system  \eqref{eq:QBODE_red} by means of Petrov-Galerkin projections. For this, we need to identify two appropriate subspaces $V\in \Rnr$ and $W\in \Rnr$ such that the state $x(t)$ can be approximated as $x(t) \approx V\hx(t)$ where $\hx(t)$ is a reduced
 state and  satisfies a Petrov-Galerkin condition as follows:
 \begin{equation*}
 W^T\left( V\dot{\hx}(t) -AV\hx(t) + HV\hx(t)\otimes V\hx(t) + \sum_{k=1}^m \hN_k \hx(t) u_k(t) + Bu(t)    \right) = 0.
 \end{equation*}
Assuming $W^TV$ is  invertible, it yields a reduced system~\eqref{eq:QBODE_red} with the realization computed as 
 \begin{equation*}
 \begin{aligned}
 \hA &= (W^TV)^{-1}W^TAV,\quad& \hN_k &= (W^TV)^{-1}W^TN_kV, ~~k \in \{1,\ldots,m\},\\
\hH &= (W^TV)^{-1}W^TH(V\otimes V), & \hB &=  (W^TV)^{-1}W^TB, \qquad \hC = CV.
 \end{aligned}
 \end{equation*}
 
 As noted in the introduction, there are several methods proposed in the literature to determine these projection matrices $V$ and $W$ or the corresponding subspaces. However, our primary objective is to extend the $\cH_2$-optimal model reduction technique \cite{morBenGG16} to the QBDAE \eqref{eq:NS_full}. Before we proceed further, we first recall important properties of the $\vecop{\cdot}$ operator, stacking the columns of a matrix on top of each other:
 \begin{equation}\label{eq:pro_vec}
\vecop{XYZ} = (Z^T\otimes X)\vecop{Y},
 \end{equation}
 where $\otimes$ denotes the Kronecker (tensor) product.  Secondly, we note down an important concept from tensor theory, that is, matricization of a tensor. 
 \begin{definition}{(e.g., \cite{koldatensor09})}
Consider a $K$-th order tensor $\cL\in \R^{I_1\times I_2\times \cdots\times I_K}$. Then, the \emph{mode-$k$ matricization} of the tensor $\cL$, denoted by $L^{(k)}$, is determined by mapping the elements  $(i_1,i_2,\ldots, i_K)$ of the tensor onto the matrix entries $(i_k,j)$ as follows:
\begin{equation*}
j  = 1+ \sum_{m=1,m\neq k}^m (i_m-1)J_m, \quad \text{where}\quad J_m := \prod_{g=1,g\neq m}^{l-1}I_g.
\end{equation*}
\end{definition}
 
Using the above concept, we can construct a 3rd-order tensor $\cH^{n\times n\times n}$ such that its mode-1 matricization of $\cH$ is the same as the Hessian $H$ in \eqref{eq:QBODE_full}, and let $H^{(2)}$ and $H^{(3)}$ denote the mode-2 and mode-3 matricization of the tensor $\cH$. Furthermore, without loss of generality, we can assume that $H(w_1\otimes w_2) = H(w_2\otimes w_1)$ for given vectors $w_1,w_2\in \Rn$, see \cite{morBenB15}. 

Having said that, we note down the definition of the $\cH_2$-norm and its truncated version, the so-called truncated $\cH_2$-norm for QB systems~\eqref{eq:QBODE_full}.  
 \begin{definition}\cite{morBenGG16}\label{def:H2_norm}
Consider a QB system \eqref{eq:QBODE_full}. Assume the $\cH_2$-norm of the system exists, the $\cH_2$-norm is defined as
\begin{equation*}
\|\Sigma\|_{\cH_2} := \sqrt{\trace{\sum_{i=1}^\infty\int_0^\infty \cdots\int_0^\infty f_i(t_1,\ldots, t_i) f_i(t_1,\ldots, t_i)^T}},
\end{equation*}
in which 
\begin{equation}
f_i(t_1,\ldots, t_i) = Cg_i(t_1,\ldots, t_i),
\end{equation}
where the $g_i(t_1,\ldots, t_i)$ satisfy
\begin{subequations}
	\begin{align*}
	g_1(t_1) & =  e^{At_1}B,\\
	g_2(t_1,t_2) & =  e^{At_2}\bbm N_1,&\ldots,& N_m\ebm\left(I_m\otimes g_1(t_1)\right),\\
	g_i(t_1,\ldots, t_i) & =  e^{At_i} \big[ H \bbm g_1(t_1)\otimes g_{i-2}(t_2,\ldots, t_{i-1}),\ldots,  g_{i-2}(t_1,\ldots, t_{i-2})\otimes g_1(t_{i-1}) \ebm,\nonumber\\
	&~~ \qquad  \bbm N_1,\ldots, N_m\ebm\left(I_m\otimes g_{i-1}(t_1,\ldots,t_{i-1})\right) \big], \qquad i\geq 3. 
	\end{align*}
\end{subequations} 
 	\end{definition}
 
 \begin{remark}
We would like to mention that the $g_i(t_1,\ldots, t_i)$ are the kernels of the Volterra series of the QB system \eqref{eq:QBODE_full} that maps the input-to-output of the system~\cite{morBenGG16}.  For more details on Volterra series expressions for QB systems, we refer to \cite{morBenGG16}. 
\end{remark}
 
 Furthermore, the connection between the above $\cH_2$-norm definition for QB systems and the recently proposed algebraic Gramians for QB systems \cite{morBTQBgoyal} has been studied in \cite{morBenGG16}. Therein, it is shown that the $\cH_2$-norm can also be computed in terms of the Gramians as follows:
 \begin{equation}
 \|\Sigma\|_{\cH_2} = \sqrt{\trace{CPC^T}} = \sqrt{\trace{B^TQB}},
 \end{equation}
where $P$ and $Q$ are the controllability and observability Gramians for QB systems which satisfy the following quadratic Lyapunov equations:
\begin{eqnarray}
AP + PA^T + H(P\otimes P)H^T + \sum_{k=1}^m N_k PN_k^T + BB^T = 0,\\
A^TQ + QA + H^{(2)}(P\otimes Q)\left(H^{(2)}\right)^T + \sum_{k=1}^m N_k^T QN_k + C^TC = 0.
\end{eqnarray} 

However,  investigating an optimal model reduction problem for QB systems using the $\cH_2$-norm is not a trivial task. Therefore, to ease the optimization problem, a concept of truncated $\cH_2$-norms for QB systems has  also been introduced in \cite{morBenGG16}, which mainly relies on the first three terms of the corresponding Volterra series. Precisely, one possible truncated $\cH_2$-norm $\|\Sigma\|_{\cH_2^{\cT}}$ is defined as follows:
\begin{equation}\label{eq:H2_trun}
\|\Sigma\|_{\cH_2^{\cT}} :=  := \sqrt{\trace{\sum_{i=1}^3\int_0^\infty \cdots\int_0^\infty \tf_i(t_1,\ldots, t_i) \tf_i(t_1,\ldots, t_i)^Tdt_1\cdots dt_i}},
\end{equation}
in which 
\begin{equation*}
\tf_i(t_1,\ldots, t_i) = C\tg_i(t_1,\ldots, t_i),
\end{equation*}
where $\tg_i(t_1,\ldots, t_i)$ satisfy
\begin{subequations}
	\begin{align*}
	\tg_1(t_1) & =  e^{At_1}B, \\
	\tg_2(t_1,t_2) & =  e^{At_2}\bbm N_1,&\ldots,& N_m\ebm\left(I_m\otimes \tg_1(t_1)\right),\\
	\tg_3(t_1,t_2,t_3) & =  e^{At_3} H \bbm \tg_1(t_1)\otimes \tg_{1}(t_2)\ebm.
	\end{align*}
\end{subequations} 

Similar to the $\cH_2$-norm expression, a connection between a truncated $\cH_2$-norm~\eqref{eq:H2_trun} and the truncated Gramians for QB systems, introduced in~\cite{morBTQBgoyal}, has also been established in \cite{morBenGG16}.  Thus, an alternative way to compute the truncated $\cH_2$-norm is as follows:
\begin{equation*}
\|\Sigma\|_{\cH_2^{\cT}}  = \sqrt{\trace{CP_\cT C^T}} = \sqrt{\trace{B^TQ_\cT Q^T}},
\end{equation*}
where $P_\cT$  are $Q_\cT$ are, respectively, truncated versions of the controllability and observability Gramians, satisfying
\begin{align*}
AP_\cT + P_\cT A^T &= -BB^T - H(P_1\otimes P_1)H^T - \sum_{k=1}^m N_kP_1N_k^T,\\
A^TQ_\cT + Q_\cT A &= -C^TC - H^{(2)}(P_1\otimes Q_1)\left(H^{(2)}\right)^T - \sum_{k=1}^m N_k^TQ_1N_k,
\end{align*}
in which $P_1$ and $Q_1$ are the unique solutions of the conventional Lyapunov equations:
\begin{align*}
AP_1 + P_1 A^T &= -BB^T,\qquad
A^TQ_1 + Q_1 A = -C^TC.
\end{align*}

Using the truncated $\cH_2$ measure, the aim is to construct a reduced-order system such that the measure of the error system is minimized. In other words, we need to determine a reduced-order system such that $\|\Sigma-\hat\Sigma\|_{\cH_2^{\cT}}$ is minimized. This problem has been studied in detail in \cite{morBenGG16}, where, first, necessary conditions for optimality are derived. Based on these conditions, an iterative scheme is proposed, which on convergence, yields a reduced-order system that \emph{approximately} satisfies  the derived optimality conditions. The iterative scheme is referred to as truncated quadratic-bilinear rational Krylov algorithm (TQB-IRKA). A brief sketch of TQB-IRKA is given in Algorithm~\ref{algo:TQBIRKA} which considers a generalized nonsingular matrix $E$ as well in the system \eqref{eq:QBODE_full} .  

\begin{algorithm}[tb!]
	\SetAlgoLined
	\DontPrintSemicolon
	\KwIn{The system matrices $E$, $A$, $H$, $N_1,\ldots,N_m$, $B$, $C$. }
	\KwOut{Reduced matrices $\hE$, $\hA$, $\hH$, $\hN_1,\ldots,\hN_m$, $\hB$, $\hC$. }
	Make an initial guess for the reduced matrices $\hE$, $\hA$, $\hH$, $\hN_1,\ldots,\hN_m$, $\hB$, $\hC$. \;
	\While{not converged}{
		Solve the generalized eigenvalue problem for the pencil $(\lambda \hE-\hA)$, i.e., determine matrices $X$ and  $Y$ such that $X\hA Y = \diag{\sigma_1,\ldots,\sigma_r} =:\Lambda$ and $X\hE Y = I_r$ in which $\sigma_i$'s are eigenvalues of the matrix pencil and $I_r$ is the identity matrix of size $r\times r$. \;
		Define $\tH = X\hH(Y\otimes Y)$, $\tN_k = X\hN_k Y$, $\tB =X\hB $ and $\tC = \hC Y$.\;
		Determine the mode-2 matricization of  $\tH$, denoted by $\tH^{(2)}$.\;
Solve for $V_1$ and $V_2$:		\label{eq:V1V2} \; 
\label{solveV1_TQB}	\quad $-EV_1\Lambda - AV_1 = B\tB^T$,\;
		\quad $-EV_2\Lambda - AV_2 = H(V_1\otimes V_1) \tH^T + \sum_{k=1}^m N_kV_1\tN_k^T$.\;
		\label{eq:W1W2}		Solve for $W_1$ and $W_2$:\;
		\quad $-E^TW_1\Lambda - A^TW_1 = C^T\tC$,\;
		\quad $-E^TW_2\Lambda - A^TW_2 = 2\cdot H^{(2)}(V_1\otimes W_1) \left(\tH^{(2)}\right)^T + \sum_{k=1}^m N_k^TW_1\tN_k$.\;
		\label{eq:VW}		Determine $V$ and $W$:\;
		\quad $V := V_1 + V_2$ and $W := W_1+W_2$.\;
		Perform $V = \ortho{V}$ and $W = \ortho{W}$.\;
		Compute reduced-order matrices:\;
		 \quad $\hE = W^TEV$,\quad  $\hA = W^TAV$,  \quad$\hN_k = W^TN_kV, ~~k \in \{1,\ldots,m\}$,\;
		 \quad $\hH = W^TH(V\otimes V)$, \quad $\hB =  W^TB$, \quad $\hC = CV.$
}
\caption{Truncated QB rational Krylov algorithm (TQB-IRKA) \cite{morBenGG16}.}
\label{algo:TQBIRKA}
\end{algorithm}

\begin{remark}
To apply Algorithm~\ref{algo:TQBIRKA}, we need to compute Kronecker products such as $H(V_1\otimes V_1)\tH^T$. In a large-scale setting, a direct computation of such Kronecker products is infeasible. As a remedy, in \cite{morBenB15,morBenGG16}, alternative methods are proposed to perform these computations  efficiently by using some tools from tensor theory or by using the special structure of the Hessian $H$, arising from semi-discretization of PDEs. 
\end{remark}

\begin{remark}
For example, Step \ref{solveV1_TQB} in Algorithm \ref{algo:TQBIRKA} which solves for $V_{1}$ is presented in a Sylvester equation form. However, using the $\vecop{\cdot}$ property \eqref{eq:pro_vec}, one can write it into a linear system, that is,
\begin{equation}
-\left( \Lambda\otimes E + I_r\otimes A\right)\vecop{V_1} = \vecop{B\tB^T},
\end{equation}
where $I_r$ is the identity matrix of size $r\times r$. Similar expressions can also be derived for the matrices $V_2$, $W_1$ and $W_2$ in Algorithm \ref{algo:TQBIRKA1}. The above $\vecop{\cdot}$ form to compute the projection matrices is very useful for the latter part of the paper. 

\end{remark}

\section{Transformation of a QBDAE into a QBODE and Model Reduction}\label{sec:transferDAEtoODE}
Our next task is to investigate how TQB-IRKA (Algorithm~\ref{algo:TQBIRKA}) can be applied to QBDAEs, having the structure as in \eqref{eq:NS_full}. For this, we first transform the system \eqref{eq:NS_full} into an equivalent ODE system by means of projections. Such a transformation is widely done in the literature for Navier-Stokes type equations, see, e.g.~\cite{MPIMD15-18,morHeiSS08}. For completeness, we show the necessary steps that transform the system \eqref{eq:NS_full} into an equivalent  ODE system. We begin by considering $B_2 = 0$ and the zero initial condition $v(0) = 0$ in \eqref{eq:NS_full}, that is:
\begin{subequations}\label{eq:NS_full_S}
	\begin{align}
	E_{11}\dot{v}(t) &= A_{11}v(t) + A_{12}p(t) + Hv(t)\otimes v(t) + \sum_{k=1}^mN_kv(t)u_k(t) + B_1u(t),\label{eq:differentialpart}\\
	0 &= A_{21}v(t),\quad v(0) = 0,\label{eq:constraint}\\
	y(t) &= C_1v(t) + C_2p(t). \label{eq:outputequation}
	\end{align}
\end{subequations}
From \eqref{eq:constraint}, we get $A_{21}\tfrac{d}{dt}v(t) = 0$, thus, leading to
\begin{align*}
0&=A_{21}\dfrac{d}{dt}v(t)  = A_{21}E_{11}^{-1}\left(E_{11}\dfrac{d}{dt}v(t)\right)\\
&=A_{21}E_{11}^{-11}\left(A_{11}v(t) + A_{12}p(t) + Hv(t)\otimes v(t) + \sum_{k=1}^mN_kv(t)u_k(t) + B_1u(t)\right)\\
&\hspace{9cm} (\text{using \eqref{eq:differentialpart}}).
\end{align*}
As a result, we obtain  an explicit expression for the pressure $p(t)$ as
\begin{equation}\label{eq:P_expression}
p(t) = -S^{-1}A_{21}E_{11}^{-1}\left(A_{11}v(t) + Hv(t)\otimes v(t) + \sum_{k=1}^mN_kv(t)u_k(t) + B_1u(t) \right),
\end{equation}
where $S = A_{21}E_{11}^{-1}A_{12}$. Substituting for $p(t)$ using~\eqref{eq:P_expression}  in \eqref{eq:differentialpart} and  \eqref{eq:outputequation} yields
\begin{subequations}\label{eq:NS_full_Pi_1}
	\begin{align}
	E_{11}\dot{v}(t) &= \Pi_l A_{11}v(t) + \Pi_l Hv(t)\otimes v(t) + \sum_{k=1}^m\Pi_l N_kv(t)u_k(t) + \Pi_l B_1u(t),\\ 
	y(t) &= \cC v(t) + \cC_H v(t)\otimes v(t) + \sum_{k=1}^m \cC_{N_k} v(t)u_k(t) + \cD u(t),\quad v(0) = 0.
	\end{align}
\end{subequations}
where 
\begin{align*}
\cC &= C_1 - C_2 S^{-1}A_{21}E_{11}^{-1}A_{11}, &\cC_H &= -C_2 S^{-1}A_{21}E_{11}^{-1}H,\\
\cC_{N_k} &=  - C_2 S^{-1}A_{21}E_{11}^{-1}N_k,  & \cD &= -C_2 S^{-1}A_{21}E_{11}^{-1}B_1,\quad \text{and}\\
\Pi_l &= I - A_{12}S^{-1}A_{21}E_{11}^{-1}. 
\end{align*}
Recall the properties of $\Pi_l$ from \cite{morHeiSS08}, which are as follows:
\begin{align*}
\Pi_l^2 &= \Pi_l, & & \range{\Pi_l} = \kernel{A_{21}E_{11}^{-1}}, & \kernel{\Pi_l}& = \range{A_{12}}.
\end{align*}
This implies that $\Pi_l$ is an \emph{oblique} projector. Moreover, for later purpose, we also define another oblique projector 
\begin{equation*}
\Pi_r = I - E_{11}^{-1}A_{12}S^{-1}A_{21}. 
\end{equation*}
First, we note down a relation between $\Pi_l$ and $\Pi_r$, that is, $\Pi_l E_{11} = E_{11}\Pi_r$. Moreover, it can be verified that 
\begin{equation*}
A_{21}z =0 \quad \text{if and only if} \quad\Pi_r z = z.
\end{equation*}
Using the properties of the projectors $\Pi_l$ and $\Pi_r$, we obtain
\begin{subequations}\label{eq:NS_full_Pi_2}
	\begin{align}
	\Pi_l E_{11}\Pi_r\dot{v}(t) &= \Pi_l A_{11}\Pi_rv(t) + \Pi_l H\left(\Pi_rv(t)\otimes \Pi_rv(t)\right) + \sum_{k=1}^m\Pi_l N_k\Pi_rv(t)u_k(t)\nonumber\hspace{-1cm}\\
	&\hspace{3cm}  + \Pi_l B_1u(t),\\ 
	y(t) &= \cC v(t) + \cC_H \Pi_rv(t)\otimes\Pi_r v(t) + \sum_{k=1}^m \cC_{N_k} \Pi_rv(t)u_k(t) + \cD u(t),
	\end{align}
\end{subequations}
with the zero initial condition, i.e., $v(0) = 0$. Moreover, the dynamics of the  system~\eqref{eq:NS_full_Pi_2} lies in the $n_v{-}n_p$ dimensional subspace, that is nothing but the null space of $\Pi_l$.  Next, we decompose the projectors $\Pi_l$ and $\Pi_r$ as follows:
\begin{equation}\label{eq:phi_theta}
\Pi_l = \theta_l\phi_l^T, \qquad \Pi_r = \theta_r\phi_r^T,
\end{equation}
in which $\theta_{\{l,r\}}$, $\phi_{\{l,r\}} \in \R^{n_v\times n_v-n_p}$ satisfy
\begin{equation*}
\theta_l^T\phi_l = I,\quad \theta_r^T\phi_r = I.
\end{equation*}
As a result, if one defines $\tv(t) := \phi_r^Tv(t)$, we consequently obtain an equivalent ODE system of the system \eqref{eq:NS_full_Pi_2}  as follows:
\begin{subequations}\label{eq:NS_full_Pi_3}
	\begin{align}
	\phi_l^T E_{11}\theta_r\dot{\tv}(t) &= \phi_l^T A_{11}\theta_r\tv(t) + \phi_l^T H\theta_r\tv(t)\otimes \theta_r\tv(t) + \sum_{k=1}^m\phi_l^T N_k\theta_r\tv(t)u_k(t)\nonumber\hspace{-1cm}\\ 
&\hspace{6cm}+ \phi_l^T B_1u(t),	\\
	y(t) &= \cC v(t) + \cC_H \theta_r\tv(t)\otimes\theta_r \tv(t) + \sum_{k=1}^m \cC_{N_k} \theta_r\tv(t)u_k(t) + \cD u(t),\label{eq:output_3}
	\end{align}
\end{subequations}
where $\tv(0)=0$.  Thus, one can apply Algorithm~\ref{algo:TQBIRKA} to obtain projection matrices which give a \emph{near}-optimal reduced system, having neglected nonlinear terms in the output equation~\eqref{eq:output_3}. However, there are two major issues: one is related to computations of $\phi_l$, $\theta_r$, or $\Pi_{\{r,l\}}$ which are expensive to compute, and in case we are able to determine these projectors and their decompositions efficiently, the matrix coefficients of the system \eqref{eq:NS_full_Pi_3}, e.g., $	\phi_l^T E_{11}\theta_r$, 	$\phi_l^T A_{11}\theta_r$, might be dense, thus making model reduction techniques numerically inefficient and expensive.  Therefore, in the following, we discuss how to employ Algorithm \ref{algo:TQBIRKA} to the system \eqref{eq:NS_full_Pi_3} without requiring  explicit computations of the projection matrices and their decompositions. 
\subsection*{Computational Issues}
In order to determine the projection matrices to compute a reduced system, we consider the following associated QB system, which is the system \eqref{eq:NS_full_Pi_3} but some neglected nonlinear terms in the output equation:
\begin{subequations}\label{eq:NS_full_Pi_4}
	\begin{align}
	\bar E\dot{\tv}(t) &= \bar A\tv(t) + \bar H \tv(t)\otimes \tv(t) + \sum_{k=1}^m\bar N_k\tv(t)u_k(t)+ \bar B_1u(t),\\ 
	\ty(t) &= \bar C \theta_r \tv(t) +  \cD u(t), \quad \tv(0) =0. 
	\end{align}
\end{subequations}
where 
\begin{equation}\label{eq:bar_matrices}
\begin{aligned}
\bar E &:= 	\phi_l^T E_{11}\theta_r,&\quad  \bar A &:= 	\phi_l^T A_{11}\theta_r , &  \quad \bar H &:= 	\phi_l^T H \theta_r \otimes \theta_r, \\
\bar N_{k} &:= 	\phi_l^T N_{k}\theta_r,&\bar B_{1} &:= 	\phi_l^T B_{1},& \bar C &:= 	 \cC \theta_r.
\end{aligned}
\end{equation}
Here, we need to neglect the nonlinear terms appearing in the output equation due to the transformation to an ODE system as Algorithm 1 is developed so far only for linear output equations. The full incorporation of the nonlinearities in the output equation requires further work. But note that the nonlinear terms are again included when the projections are applied to obtain the reduced-order model, see Remark \ref{remark:proj_nonlinear}.

As a first step, we aim at determining the projection matrices $V$ and $W$ such that the original matrices like $E_{11}$, $A_{11}$ can be used to compute the reduced matrices instead of using, e.g., $\bar E$, $\bar A$, i.e., a reduced matrix $\hE$ can be computed as  $W^T E_{11}V$, and so on.  

For this purpose, let $\bar V$ and $\bar W$ be the solutions of the Sylvester equations in Algorithm~\ref{algo:TQBIRKA} (steps 6--\ref{eq:VW}) using the matrices $\bar E$, $\bar A$, etc.  Furthermore, we define the matrices $V$ and $W$, satisfying:

\begin{equation}
V = \theta_r\bar V, \qquad \text{and}\qquad W = \phi_l\bar W.
\end{equation}
Then, it can  easily be verified  that the reduced matrices computed using the quantities denoted by \emph{bar}, e.g., $\bar V$, $\bar W$, $\bar E$ are the same as the reduced matrices computed using $V$ and $W$ and the original matrices. In other words, $\bar W^T \bar E \bar V = W^T E_{11}V$ and so on. Using this formulation, we sketch  an algorithm for model reduction of \eqref{eq:NS_full_Pi_4} which gives near-optimal reduced systems on convergence by projecting the original system matrices, see Algorithm \ref{algo:TQBIRKA1}.

\begin{algorithm}[tb!]
	\SetAlgoLined
	\DontPrintSemicolon
	\KwIn{The system matrices $E_{11}$, $A_{11}$, $H$, $n_v,\ldots,N_m$, $B_1$, $\cC$. }
	\KwOut{Redcued matrices $\hE$, $\hA$, $\hH$, $\hn_v,\ldots,\hN_m$, $\hB$, $\hC$. }
	Make an initial guess of reduced matrices $\hE$, $\hA$, $\hH$, $\hn_v,\ldots,\hN_m$, $\hB$, $\hC$. \;
	\While{not converged}{
		Solve the generalized eigenvalue problem for the pencil $(\lambda \hE-\hA)$, i.e., determine matrices $X$ and  $Y$ such that $X\hA Y = \diag{\sigma_1,\ldots,\sigma_r} =:\Lambda$ and $X\hE Y = I_r$ in which the $\sigma_i$'s are eigenvalues of the matrix pencil.  \;
		Define $\tH = X\hH(Y\otimes Y)$, $\tN_k = X\hN_k Y$, $\tB =X\hB $ and $\tC = \hC Y$.\;
		Determine  themode-2 matricization of  $\tH$, denoted by $\tH^{(2)}$.\;
		Define $L := \left(-\Lambda \otimes \bar E - I_r \otimes \bar A\right)^{-1}$, where $\bar E$ and $\bar A$ are as in \eqref{eq:bar_matrices}.\;
Solve for $V_1$ and $V_2$:\; 
	\label{eq:V1V2_1}	\quad $\vecop{V_1} = (I_r\otimes \theta_r) L(I_r\otimes \phi_l^T)\vecop{B\tB^T}$, \;
		\quad $\vecop{V_2} = (I_r\otimes \theta_r) L(I_r\otimes \phi_l^T)\vecop{H(V_1\otimes V_1) \tH^T + \sum_{k=1}^m N_kV_1\tN_k^T}$. \;
	Solve for $W_1$ and $W_2$:\;
	\label{eq:W1W2_1}	\quad $\vecop{W_1} = (I_r\otimes \phi_l) L^T(I_r\otimes \theta_r^T)\vecop{\cC^T\tC}$, \;
		\quad $\vecop{W_2} = (I_r\otimes \phi_l) L^T(I_r\otimes \theta_r^T)\vecop{2\cdot H^{(2)}(V_1\otimes W_1) \left(\tH^{(2)}\right)^T + \sum_{k=1}^m N_k^TW_1\tN_k}$. \;
	\label{eq:VW_1}	Determine $V$ and $W$:\;
		\quad $V := V_1 + V_2$ and $W := W_1+W_2$.\;
		Perform $V = \ortho{V}$ and $W = \ortho{W}$.\;
		Compute reduced-order matrices:\;
		\quad $\hE = W^TE_{11}V$,\quad  $\hA = W^TA_{11}V$,  \quad$\hN_k = W^TN_kV, ~~k \in \{1,\ldots,m\}$,\;
		\quad $\hH = W^TH(V\otimes V)$, \quad $\hB =  W^TB_1$, \quad $\hC = \cC V.$
	}
	\caption{Truncated QB rational Krylov algorithm for the system~\eqref{eq:NS_full_Pi_4}, involving projectors.}
	\label{algo:TQBIRKA1}
\end{algorithm}


However, to compute the projection matrices at each iteration in Algorithm \ref{algo:TQBIRKA1}, we still require the basis matrices  $\phi_l$ and $\theta_r$. Therefore, our next goal is to determine these projection matrices without involving $\phi_l$ and $\theta_r$, or more precisely, we aim at computing the matrices $V$ and $W$ using only the original matrices such as $E_{11}$, $A_{11}$, $A_{12}$. 

Luckily, a similar problem has been studied in \cite{goyal2016iterative}, where for a symmetric case, i.e., $\Pi_l = \Pi_r^T$, it is shown  how to compute $(I_r\otimes \theta_r) L(I_r\otimes \phi_l^T)f$ for  given arbitrary vector $f$ and matrix $L$, without explicitly forming $\theta_r$ and $\phi_l$. However, for the unsymmetric case, i.e., $\Pi_l \neq \Pi_r^T$, the following result can be developed in a similar fashion as  in \cite{goyal2016iterative}, and, therefore,  a detailed proof is omitted. 

\begin{lemma}\label{eq:lemma_saddlepoint}
Consider $\phi_l$ and $\theta_r$ as in \eqref{eq:phi_theta}, and assume  $\cX = (I_r\otimes \phi_l^T) \cT (I_r\otimes \theta_r)$ is invertible for a given $\cT$. Furthermore, let $\cG$ and $\cG_T$  be  $ (I_r\otimes \theta_r) \cX^{-1}(I_r\otimes \phi_l^T)$ and $ (I_r\otimes \theta_l) \cX^{-T}(I_r\otimes \phi_r^T)$, respectively. Then, the vector 
\begin{equation*}
\bar v = \cG f
\end{equation*} 
solves
\begin{equation*}
\bbm \cT & I_r\otimes A_{12} \\ I_r\otimes A_{21} & 0 \ebm \bbm \bar v \\ \xi_v \ebm= \bbm f \\ 0 \ebm. 
\end{equation*}
Similarly, the vector 
\begin{equation*}
\bar w = \cG_T g
\end{equation*} 
solves
\begin{equation*}
\bbm \cT^T & I_r\otimes A_{21}^T \\ I_r\otimes A_{12}^T & 0 \ebm \bbm \bar w \\ \xi_w \ebm = \bbm g \\ 0 \ebm. 
\end{equation*}
\end{lemma}

By making use of the result of Lemma \ref{eq:lemma_saddlepoint}, it can be shown that the steps 7-10 in Algorithm \ref{algo:TQBIRKA1} can be performed without an explicit requirement of the basis matrices $\phi_l$ and $\theta_r$. We rather need to solve appropriate saddle point problems to compute $V_{\{1,2\}}$ and $W_{\{1,2\}}$. All these analyses lead to Algorithm~\ref{algo:TQBIRKA2} that does not require any explicit computation of the basis matrices $\phi_l$ and $\theta_r$ even for computing projection matrices. 

\begin{algorithm}[tb!]
	\SetAlgoLined
	\DontPrintSemicolon
	\KwIn{The system matrices $E_{11}$, $A_{11}$, $H$, $n_v,\ldots,N_m$, $B_1$, $\cC$. }
	\KwOut{Redcued matrices $\hE$, $\hA$, $\hH$, $\hn_v,\ldots,\hN_m$, $\hB$, $\hC$. }
	Make an initial guess of reduced matrices $\hE$, $\hA$, $\hH$, $\hn_v,\ldots,\hN_m$, $\hB$, $\hC$. \;
	\While{not converged}{
		Solve the generalized eigenvalue problem for the pencil $(\lambda \hE-\hA)$, i.e., determine matrices $X$ and  $Y$ such that $X\hA Y = \diag{\sigma_1,\ldots,\sigma_r} =:\Lambda$ and $X\hE Y = I_r$ in which $\sigma_i$'s are eigenvalues of the matrix pencil.  \;
		Define $\tH = X\hH(Y\otimes Y)$, $\tN_k = X\hN_k Y$, $\tB =X\hB $ and $\tC = \hC Y$.\;
		Determine mode-2 matricization of  $\tH$, denoted by $\tH^{(2)}$.\;
		Define $L := \bbm \left(-\Lambda\otimes E_{11} - I_r\otimes A_{11}  \right)^{-1} & I_r\otimes A_{12} \\ I_r\otimes A_{21} & 0 \ebm $, and $\cI = \bbm I_{n_v} & 0 \ebm$.\;
		Solve for $V_1$ and $V_2$:\; 
		\label{eq:V1V2_1}	\quad $\vecop{V_1} = \cI L\bbm \vecop{B\tB^T}\\ 0 \ebm$, \;
		\quad $\vecop{V_2} = \cI L \bbm \vecop{H(V_1\otimes V_1) \tH^T + \sum_{k=1}^m N_kV_1\tN_k^T}\\ 0 \ebm$. \;
		Solve for $W_1$ and $W_2$:\;
		\label{eq:W1W2_1}	\quad $\vecop{W_1} = \cI L^T \bbm \vecop{\cC^T\tC}\\ 0 \ebm$, \;
		\quad $\vecop{W_2} = \cI L^T \bbm \vecop{2\cdot H^{(2)}(V_1\otimes W_1) \left(\tH^{(2)}\right)^T + \sum_{k=1}^m N_k^TW_1\tN_k} \\ 0 \ebm$. \;
		\label{eq:VW_1}	Determine $V$ and $W$:\;
		\quad $V := V_1 + V_2$ and $W := W_1+W_2$.\;
		Perform $V = \ortho{V}$ and $W = \ortho{W}$.\;
		Compute reduced-order matrices:\;
		\quad $\hE = W^TE_{11}V$,\quad  $\hA = W^TA_{11}V$,  \quad$\hN_k = W^TN_kV, ~~k \in \{1,\ldots,m\}$,\;
		\quad $\hH = W^TH(V\otimes V)$, \quad $\hB =  W^TB_1$, \quad $\hC = \cC V.$
	}
	\caption{Truncated QB rational Krylov algorithm for the system~\eqref{eq:NS_full_Pi_4}.}
	\label{algo:TQBIRKA2}
\end{algorithm}
 \begin{remark}\label{remark:proj_nonlinear}
Recall that Algorithm \ref{algo:TQBIRKA2} gives on convergence a \emph{near}-optimal reduced system for the system \eqref{eq:NS_full_Pi_3}, having neglected the nonlinear terms in the output equation in the model reduction process. Nonetheless, we reduce these nonlinear terms in the output equation using the projection matrix $V$, obtained on convergence, i.e., $\hat{\cC}_H = \cC_H(V\otimes V)$ and $\hat{\cC}_{N_k} = \cC_{N_k}V$. Furthermore, if the output of the system \eqref{eq:NS_full_S} is given only by  linear combinations of the velocity $v(t)$, then all nonlinear terms in the output equation of the system~\eqref{eq:NS_full_Pi_3} are zero.  
 	\end{remark}
\begin{remark}
Throughout the above discussion, we have assumed that $B_2 =0$ in \eqref{eq:NS_full}. However, as discussed, e.g., in \cite{morHeiSS08}, the general case $B_2\neq 0$ can be converted into the case $B_2= 0$  by an appropriate change of variables for $v(t)$. For this, one needs to decompose $v(t)$ as follows:
\begin{equation}\label{eq:decompose_v}
v(t) = v_m(t) + v_u(t),
\end{equation}
where $v_u(t) =  \Omega u(t)$ in which $\Omega := -E_{11}^{-1}A_{12}(A_{21}E_{11}A_{12})B_2$.  If one substitutes for $v(t)$ from \eqref{eq:decompose_v} into \eqref{eq:NS_full}, then it can be easily seen that $A_{21}v_m(t) =0$. Furthermore, we assume  an initial condition to be consistent. Performing  similar algebraic calculations as done for the problem $B_2 =0$ leads to the following system:
\begin{subequations}\label{eq:NS_full_Gen}
	\begin{align}[]
	\Pi_l E_{11}\Pi_r\dot{v}(t) &= \Pi_l A_{11}\Pi_rv(t) + \Pi_l H\left(\Pi_rv(t)\otimes \Pi_rv(t)\right) + \sum_{k=1}^m\Pi_l \cN_k\Pi_rv(t)u_k(t)\nonumber\hspace{-1cm}\\
	&\hspace{3cm}  + \Pi_l\cB_1u(t) + \Pi_l \cB_u (u(t)\otimes u(t)),\\ 
	y(t) &= \cC v(t) + \cC_H \Pi_rv(t)\otimes\Pi_r v(t) + \sum_{k=1}^m \cC_{N_k} \Pi_rv(t)u_k(t) + \cD u(t) \nonumber\hspace{-1cm}\\
	& \hspace{3cm} - C_2 (A_{21}E_{11}^{-1}B_2)\dot u(t),
	\end{align}
\end{subequations}
where 
\begin{align*}
\cN_k &= N_k - H(I\otimes \Omega_k +  \Omega_k\otimes I) , & \cB_1 &= B_1 + A_{11}\Omega,\\ \cB_u&= H(\Omega\otimes \Omega) +  \bbm n_v\Omega, \ldots, N_m\Omega\ebm, & \cC& = C_1 - C_2 S^{-1}A_{21}E_{11}^{-1}A_{11}, \\
\cC_H  &= -C_2 S^{-1}A_{21}E_{11}^{-1}H, & \cC_{N_k} &=  -C_2 S^{-1}A_{21}E_{11}^{-1}\cN_k,
\end{align*}
in which $S := A_{21}E_{11}^{-1}A_{12}$. 
Thus, we obtain a  system equivalent to \eqref{eq:NS_full} which has a similar form as in \eqref{eq:NS_full_Pi_2}. However, the system \eqref{eq:NS_full_Gen} contains some extra terms which are functions of the input $u(t)$, e.g., $u_p\cdot u_q$, $\{p,q\}\in \{1,\ldots,m\}$ and the derivative of the input $u(t)$. Although they are  functions of the input in the forward simulation, we consider them formally as different inputs as far as the model reduction problem is concerned. Hence, Algorithm \ref{algo:TQBIRKA2} can readily be applied to the system \eqref{eq:NS_full_Gen} to determine reduced-order systems. 

\end{remark}

\section{Numerical Experiment}\label{sec:numerics}
In this section, we test the efficiency of Algorithm \ref{algo:TQBIRKA2} using a numerical example of a lid-driven cavity, obtained using semi-discretization of the Navier--Stokes equation.  We initialize Algorithm \ref{algo:TQBIRKA2} randomly and iterate it until the relative change in the eigenvalues of the reduced matrix $\hA$ is less than $10^{-5}$. All the simulations were done on a board with 4 \intel~\xeon ~E7-8837 CPUs with a 2.67-
GHz clock speed using \matlab~8.0.0.783 (R2012b). 

\subsection{Lid Driven Cavity}
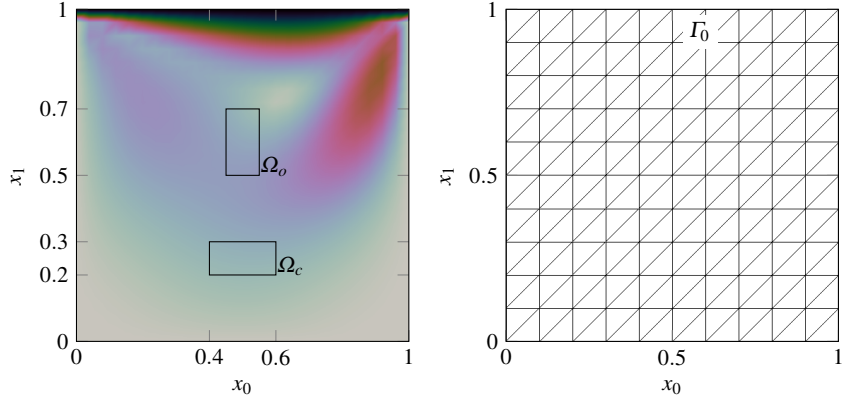
\begin{figure}[tb]
	\newlength\figureheight
	\newlength\figurewidth
	\setlength\figureheight{6cm}
	\setlength\figurewidth{6cm}
%
%
%
%
\begin{tikzpicture}

\begin{axis}[
xmin=0, xmax=1,
ymin=0, ymax=1,
width=\figureheight,
height=\figurewidth,
axis on top,
xtick={0,0.4,0.6,1},
ytick={0,.2,.3,.5,.7,1},
tick label style={font=\footnotesize},
xlabel={$x_0$},
ylabel={$x_1$},
xlabel near ticks,
ylabel near ticks
]
\addplot graphics [includegraphics cmd=\pgfimage,xmin=0, xmax=1, ymin=0, ymax=1] {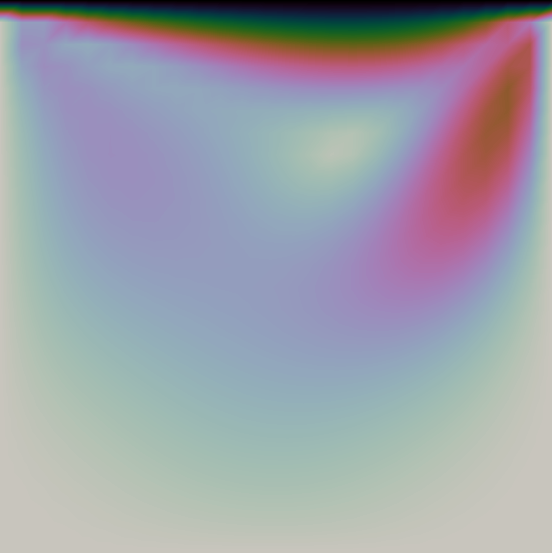};
\path [draw=black, fill opacity=0] (axis cs:13,1)--(axis cs:13,1);

\path [draw=black, fill opacity=0] (axis cs:0.05,13)--(axis cs:0.05,13);

\path [draw=black, fill opacity=0] (axis cs:13,1.38777878078145e-17)--(axis cs:13,1.38777878078145e-17);

\path [draw=black, fill opacity=0] (axis cs:0,13)--(axis cs:0,13);

\draw [color=black] (axis cs:0.4,0.2) rectangle (axis cs:0.6,0.3);
\node at (axis cs:0.58,0.28) [color=black, anchor=north west] {$\Omega_c$};
\draw [color=black] (axis cs:0.45,0.5) rectangle (axis cs:0.55,0.7);
\node at (axis cs:0.53,0.58) [color=black, anchor=north west] {$\Omega_o$};
\end{axis}

\end{tikzpicture}
%
%
%
%
\begin{tikzpicture}

\begin{axis}[
xmin=0, xmax=1,
ymin=0, ymax=1,
width=\figureheight,
height=\figurewidth,
axis on top,
xtick={0,.5,1},
ytick={0,.5,1},
tick label style={font=\footnotesize},
xlabel={$x_0$},
ylabel={$x_1$},
xlabel near ticks,
ylabel near ticks
]
\addplot graphics [includegraphics cmd=\pgfimage,xmin=0, xmax=1, ymin=0, ymax=1] {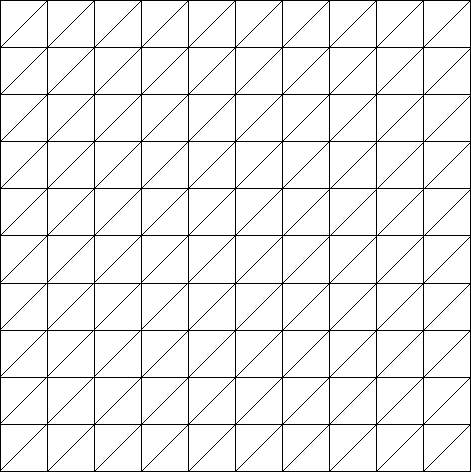};

\path [draw=black, fill opacity=0] (axis cs:13,1)--(axis cs:13,1);

\path [draw=black, fill opacity=0] (axis cs:0.05,13)--(axis cs:0.05,13);

\path [draw=black, fill opacity=0] (axis cs:13,1.38777878078145e-17)--(axis cs:13,1.38777878078145e-17);

\path [draw=black, fill opacity=0] (axis cs:0,13)--(axis cs:0,13);

\node at (axis cs:0.53,0.88) [color=black, fill=white, anchor=south west] {$\Gamma_0$};

\end{axis}

\end{tikzpicture}
	\caption{Illustration of the velocity magnitude of the lid-driven cavity
		problem and the domain of control and observation
		$\Omega_c=[0.4,0.6]\times [ 0.2,0.3]$ and $\Omega_o = [0.45,0.55] \times
		[0.5,0.7]$.}
	\label{fig:drivcav}
\end{figure}

Here, we consider a lid-driven cavity as shown in Figure \ref{fig:drivcav}, which is modeled using the incompressible Navier-Stokes equations in the velocity $\breve{v}$ and the pressure $\breve{p}$ on the unit square $\Omega = (0,1)^2$. The governing equations and boundary conditions are:
\begin{subequations}
\begin{align*}
\dot{\breve{v}} + \left(\breve{v}\cdot \nabla \right) \breve{v} - \dfrac{1}{\Rey} \Delta \breve{v} +  \nabla \breve{p}&=0,\\
\nabla \cdot \breve{v} &= 0,
\end{align*}
\end{subequations}
and boundary and initial conditions are as follows:
\begin{equation*}
\left.\breve{v}\right|_{\text{at boundary}} = g, \qquad \left.\breve{v}\right|_{t = 0} = \breve{v}_0.,
\end{equation*}
where $\Rey$ is the Reynolds number, $g$ represents the Dirichlet boundary conditions, and we set $\breve{v} = \bbm 1 \\ 0 \ebm $ for the upper boundary and zero at the other boundary, and $\breve{v}_0$ is an initial condition. Furthermore, the system is subject to a control input $u(t)$ in the domain $\Omega_c$. As a  result, we ensure that the $x$  and $y$ components of the  velocity in the domain $\Omega_c$ are $u(t)$. 

\begin{figure}[tb]
	\centering
	\includegraphics[trim={0cm 4.1cm 1cm 4.2cm},clip, height = 4.5cm, width = 4.5cm]{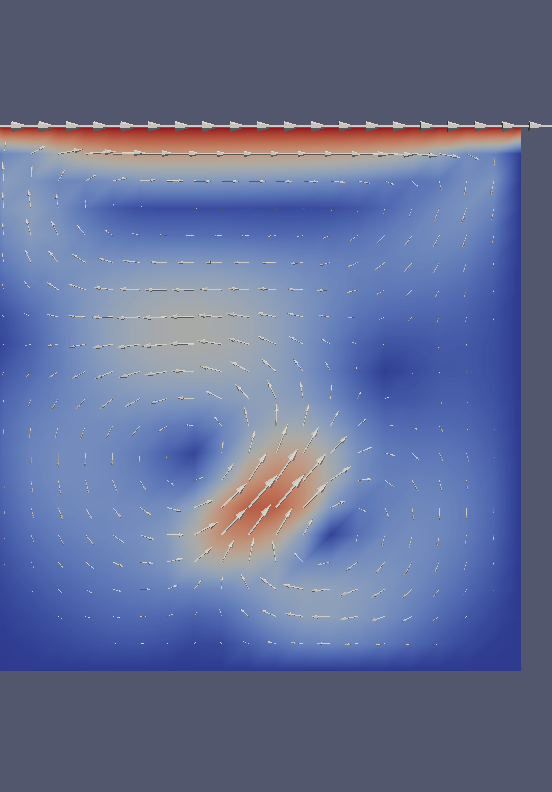}
		\includegraphics[trim={0 1.4cm 0cm 1.4cm},clip,height = 4.5cm, width = 4.5cm]{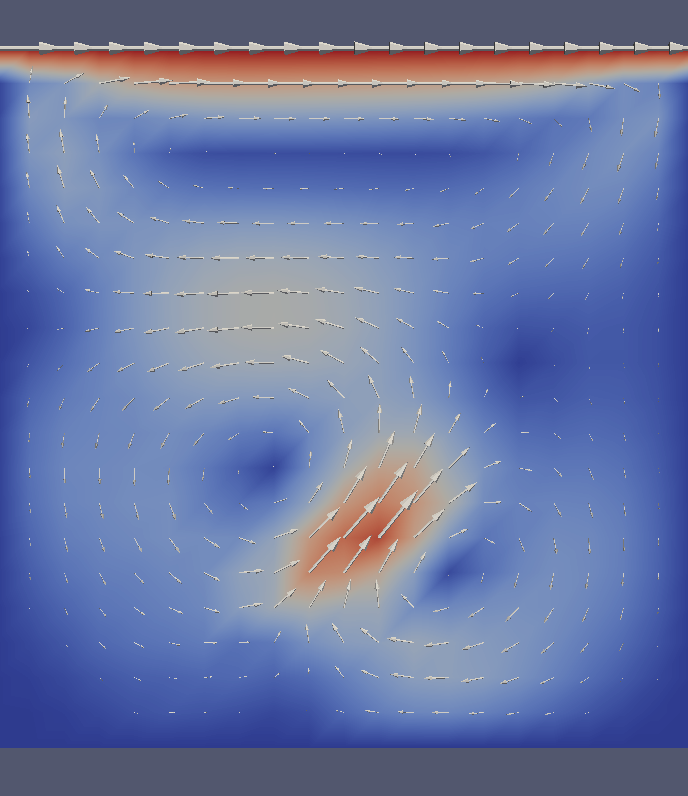}
		\caption{Comparison of $|v|$ obtained from the original (left) and reduced (right) systems at time $t  =3.75$s for an input $u(t) = 2t^2\exp(-t/2)\sin(2\pi t/5)$.}
		\label{fig:fullgrid}
\end{figure}

A similar problem has been considered in the model reduction framework in \cite{MPIMD15-18} for a SISO system and with a different control setting. Having applied a finite element scheme using the Taylor-Hood scheme on a uniform mesh, we obtain a discretized system in velocity $v$ and pressure $p$ as follows:
\begin{subequations}
\begin{align*}
E_{11}\dot{v} &= A_{11}v(t) +  H(v\otimes v) +  A_{12}v(t) + f + B_1 u(t),\\
0& = A_{12}^T v(t),\qquad v(0) = v_s,
\end{align*}
\end{subequations}
where $E_{11}, A_{11} \in \R^{n_v\times n_v} $, $A_{12} \in \R^{n_v\times n_p}$, $H \in \R^{n_v\times n_v^2}$, $f \in \R^{n_v}$ and $v_s$ is the steady-state solution of the system when the control input $u$ is set to zero. Furthermore, we observe the velocity in the domain $\Omega_o$, leading to an output equation as follow:
\begin{equation*}
y(t) = C v(t).
\end{equation*}
Since the initial condition is considered to be the steady-state solution of the system, we perform a change of variables as $v_\delta =  v +  v_s$ and $p_\delta = p + p_s$ to ensure the zero initial condition of the transformed system. This results in the following system:
\begin{subequations}
	\begin{align*}
	E_{11}\dot{v}_\delta &= \left(A_{11}+ X\right)v_\delta(t) +  H(v_\delta\otimes v_\delta) +  A_{12}p_\delta(t)  + B_1 u(t),\\
	0& = A_{12}^T v_\delta(t),\qquad v_\delta(0) = 0,
	\end{align*}
\end{subequations}
where $X :=  H(v_s\otimes I + I\otimes v_s)$. We choose the degrees of freedom for the velocity and the pressure $3042$ and $440$, respectively, i.e.,  $n_v = 3042$ and $n_p = 440 $. We set the Reynolds number $\Rey=100$.  Furthermore, we select four nodes in the domain $\Omega_o$ on which we measure the $x$ and $y$ components of the velocity; thus, we have $8$ outputs.  Next, we employ Algorithm \ref{algo:TQBIRKA2}  to obtain a reduced-order system of order $r  =140$. To check the accuracy of the reduced system, we perform time-domain simulations for the original and reduced systems for an arbitrary control input $u(t) = 2t^2\exp(-t/2)\sin(2\pi t/5)$ and observe the outputs which are plotted in Figure \ref{fig:response}. This shows that the reduced system captures the dynamics of the system faithfully. Furthermore, we also plot the velocity at the full grid in Figure~\ref{fig:fullgrid}, although the considered model reduction aims at capturing only the input-to-output behavior of the system. The figure indicates that the reduced-order system can replicate the visual dynamics of the original system at the full grid as well. 
\newlength{\fheight}
\newlength{\fwidth}

\definecolor{mycolor1}{rgb}{0.00000,0.75000,0.75000}%
\definecolor{mycolor2}{rgb}{0.75000,0.00000,0.75000}%
\definecolor{mycolor3}{rgb}{0.75000,0.75000,0.00000}%
\newcommand{\hwplotA}{\raisebox{3pt}{\tikz{\draw[color=black!50!green,line width=1.2pt](0,0) -- (5mm,0);}} \raisebox{3pt}{\tikz{\draw[color=red,line width=1.2pt](0,0) -- (5mm,0);}} \raisebox{3pt}{\tikz{\draw[color=mycolor1,line width=1.2pt](0,0) -- (5mm,0);}} \raisebox{3pt}{\tikz{\draw[color=mycolor2,line width=1.2pt](0,0) -- (5mm,0);}}  \raisebox{3pt}{\tikz{\draw[color=mycolor3,line width=1.2pt](0,0) -- (5mm,0);}} \raisebox{3pt}{\tikz{\draw[color=darkgray,line width=1.2pt](0,0) -- (5mm,0);}}  \raisebox{3pt}{\tikz{\draw[color=blue,line width=1.2pt](0,0) -- (5mm,0);}} }

\newcommand{\defsp}{\hspace{.38cm}}
\newcommand{\hwplotB}{\defsp\raisebox{1pt}{\tikz{\draw[line width = 1.2pt,color =black!50!green ] (2.5mm,0) circle (2pt);}} \defsp\raisebox{1pt}{\tikz{\draw[color=red,line width=1.2pt] (2.5mm,0) circle (2pt);}}\defsp\raisebox{1pt}{\tikz{\draw[color=mycolor1,line width=1.2pt] (2.5mm,0) circle (2pt);}} \defsp\raisebox{1pt}{\tikz{\draw[color=mycolor2,line width=1.2pt] (2.5mm,0) circle (2pt);}} \defsp \raisebox{1pt}{\tikz{\draw[color=mycolor3,line width=1.2pt] (2.5mm,0) circle (2pt);}} \defsp\raisebox{1pt}{\tikz{\draw[color=darkgray,line width=1.2pt] (2.5mm,0) circle (2pt);}} \defsp \raisebox{1pt}{\tikz{\draw[color=blue,line width=1.2pt] (2.5mm,0) circle (2pt);}} }

\begin{figure}[tb]
	\centering
	\begin{tikzpicture}
	\begin{customlegend}[legend columns=1, legend style={/tikz/every even column/.append style={column sep=0.5cm}} , legend entries={ \hwplotA Original system $(n = 3482)$, \hspace{-1.3cm}\hwplotB~~ \emph{TQB-IRKA} $(r = 140)$,~BT $(r = 20)$}, ]
	\addlegendimage{blue,line width = 1.2pt}     \addlegendimage{blue, only marks,mark = o,line width = 1.2pt}
	\end{customlegend}
	\end{tikzpicture}
		\setlength\fheight{4.0cm}
		\setlength\fwidth{7.25cm}
	\input{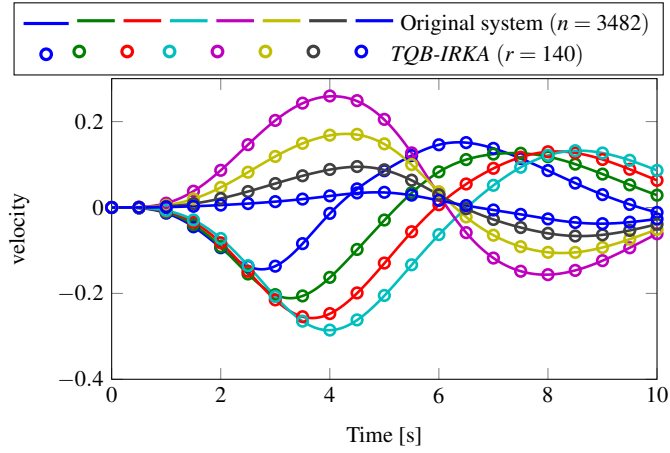}
	 \caption{For a random control $u(t) = 2t^2\exp(-t/2)\sin(2\pi t/5)$, we compare the $x$ and $y$ components of velocities at four nodes in the domain $\Omega_o$  obtained via the original and reduced systems in the figure which are indicated by different colors. }
	 \label{fig:response}
\end{figure}

\section{Conclusions}\label{sec:conclusions}

In this work, we have studied a model reduction problem for a particular class of quadratic-bilinear descriptor systems, especially arising from semi-discretization of the Navier-Stokes equations. In particular, we have investigated how one can employ the iterative model reduction scheme proposed in \cite{morBenGG16} for quadratic-bilinear ODEs to the latter class of descriptor systems. To that end, we have first transformed the quadratic-bilinear descriptor system into an equivalent ODE system by means of projectors. Furthermore,  for practical computations, we have proposed an efficient iterative scheme for the considered quadratic-bilinear descriptor systems to avoid the explicit computation of the projectors that are used to transform the system into an ODE system. Using as numerical example a semi-discretized Navier-Stokes equations, we have shown the efficiency of the proposed method. 

As a further research topic, an extension of the balanced truncation method for quadratic-bilinear ODE systems \cite{morBTQBgoyal} to  descriptor systems will be very useful. Furthermore, 
it will be a significant contribution to extended these methods to other classes of descriptor systems, particularly, appearing in mechanical systems. 

 \bibliographystyle{siam}
  \bibliography{mor.bib}
 
 \end{document}